\documentclass[12pt,leqno]{amsart}
\usepackage{amssymb}
\usepackage{amsmath,amssymb,color}
\numberwithin{equation}{section}

\newcommand{\la}{\lambda}
\newcommand{\va}{\varphi}
\newcommand{\ppp}{\partial}

\newcommand{\hhalf}{\frac{1}{2}}

\newcommand{\ddd}{\mbox{div}\thinspace}

\newcommand{\wwwu}{\widetilde{u}}
\newcommand{\wwwv}{\widetilde{v}}

\newcommand{\weight}{e^{2s\va}}

\newcommand{\R}{\mathbb{R}}

\newcommand{\www}{\widetilde}

\newcommand{\ooo}{\overline}
\newcommand{\OOO}{\Omega}

\newcommand{\BBB}{\mathcal{B}}

\allowdisplaybreaks

\title
[]
{
Unique continuation for a mean field game system}

\author{$^1$ Oleg Imanuvilov,  
$^2$ Hongyu Liu and $^3$ Masahiro Yamamoto}

\thanks{
$^1$ Department of Mathematics, Colorado State University, 101 Weber Building, 
Fort Collins CO 80523-1874, USA 
e-mail: {\tt oleg@math.colostate.edu}
\\
$^2$ Department of Mathematics, City University of Hong Kong, Kowloon, 
Hong Kong SAR, China email: {\tt hongyliu@cityu.edu.hk}
\\
$^3$ Graduate School of Mathematical Sciences, The University
of Tokyo, Komaba, Meguro, Tokyo 153-8914, Japan 
e-mail: {\tt myama@ms.u-tokyo.ac.jp}
}

\date{}
\begin{document}
\maketitle

\begin{abstract}
For a mean field game system, 
we prove the unique continuation which asserts that if Cauchy 
data are zero on arbitrarily chosen lateral subboundary, then 
the solution identically vanishes.
\end{abstract} 
\baselineskip 18pt

\section{Introduction and key Carleman estimate}

Let $\OOO \subset \R^n$ be a bounded domain with smooth boundary 
$\ppp\OOO$, let $T>0$, and $Q:= \OOO\times (0,T)$.
In this article, we consider a system of
the mean field game:
$$
\left\{ \begin{array}{rl}
& \ppp_tu(x,t) + a_1(x,t)\Delta u(x,t)
- \hhalf\kappa(x,t)\vert \nabla u(x,t)\vert^2
- h(x,t)u = F, \\
& \ppp_tv(x,t) - \Delta (a_2(x,t)v) 
- \ddd (\kappa(x,t)v(x,t)\nabla u(x,t)) = G \quad 
\mbox{in $Q$}.
\end{array}\right.
                                       \eqno{(1.1)}
$$
Throughout this article, we assume 
$$
a_1, a_2 \in C^{2}(\ooo{Q}), \, > 0 \quad \mbox{on $\ooo{Q}$},
\quad \kappa \in C^{1,0}(\ooo{Q}), \quad h\in L^{\infty}(Q).
                                   \eqno{(1.2)}
$$
Here and henceforth we set 
$C^{2,2}(\ooo{Q}):= \{ a\in C(\ooo{Q});\, \ppp_t^ka, \ppp_i\ppp_ja,
\ppp_ia, a \in C(\ooo{Q})$ for $0\le k \le 2$ \\
and $1\le i,j\le n \}$, 
$C^{1,0}(\ooo{Q}):= \{ a\in C(\ooo{Q});\, \partial_j a \in C(\ooo{Q}),\quad 1\le j\le n\}$
and \\
$H^{2,1}(Q):= \{ w \in L^2(Q);\, 
w, \ppp_tw, \ppp_iw, \ppp_i\ppp_jw \in L^2(Q) \,\,
\mbox{for $1\le i,j \le n$}\}$.

In (1.1), $x$ and $t$ are the state and the time variables, 
and $u$ and $v$ denote the value of the game and the population 
density of players respectively (e.g., Achdou, Cardaliaguet, Delarue,
Porretta and Santambrogio \cite{ACDPS}, Lasry and Lions \cite{LL}).

The main purpose of this article is establish the unique continuation for
(1.1):
\\
{\bf Theorem 1.} 
{\it 
Let $\gamma \subset \ppp\OOO$ be arbitrarily chosen
non-empty relatively open subboundary.
We assume that $(u,v), (\wwwu, \wwwv) \in H^{2,1}(Q)\times H^{2,1}(Q)$ 
satisfy (1.1) and
$$\left\{\begin{array}{rl}
& u, \nabla u, \Delta u, \wwwu, \nabla\wwwu, \Delta \wwwu \in L^{\infty}(Q),\\
& v, \nabla v, \wwwv, \nabla\wwwv \in L^{\infty}(Q), \quad 
\ppp_t(u-\wwwu), \ppp_t(v-\wwwv) \in L^2((\ppp\OOO\setminus \Gamma) 
\times (0,T)).
\end{array}\right.
                                       \eqno{(1.3)}
$$
Then $u=\wwwu$, $\nabla u = \nabla \wwwu$, $v=\wwwv$ and $\nabla v
= \nabla \wwwv$ on $\Gamma \times (0,T)$ implies $u=\wwwu$ and $v=\wwwv$ 
in $Q$.
}

The unique continuation for a single parabolic equation is well-known and as
early works we can refer to Mizohata \cite{Miz}, Saut and Scheurer \cite{SS}
among many other works.

The mean field game system is a mixed type of backward and forward parabolic
equations for $u$ and $v$ respectively, and so the initial boundary value 
problem requires special cares.  However, it is straightforward to 
establish a Carleman estimate which is applicable to the unique continuation 
and other problems such as inverse problems.

It is known that a relevant Carleman estimate can derive the conditional 
stability in determining $u$ and $v$ by Cauchy data on $\Gamma \times (0,T)$.
For example we can refer to Klibanov and Timonov \cite{KT} under
geometric constraints on $\Gamma$.  We can establish the conditional 
stability from arbitrary subboundary by the same way as 
Huang, Imanuvilov and Yamamoto \cite{HIY}, but here we omit the details.
  
Our key Carleman estimate is stated as Theorem 2 in Section 2, and 
is derived directly thanks to that the second-order coupling terms of 
$u$ appear in the equation in $v$ as $\Delta u$.  For general cases,
such coupling should be a linear combination of $\ppp_i\ppp_ju$, $1\le i,j\le n$.
The derivation of a relevant Carleman estimate is more complicated and 
in a forthcoming work we will pursue.  

As for inverse problems, we refer to Klibanov
\cite{Kl23}, Klibanov and Averboukh \cite{KlAv},
Klibanov, Li and Liu \cite{KLL1}, \cite{KLL2},
Liu, Mou and Zhang \cite{LMZ}, Liu and Yamamoto \cite{LY},
Liu and Zhang \cite{LZ1}, \cite{LZ2} 

This article is composed of three sections.  In Section 2, we prove a key 
Carleman estimate for (1.1) and Section 3 is devoted to the completion of 
the proof of Theorem 1.  

\section{Key Carleman estimate}

For subboundary $\Gamma \subset \ppp\OOO$, we see that there exists
$d \in C^2(\ooo{\OOO})$ such that 
$$
d>0 \quad \mbox{in $\OOO$}, \quad \vert \nabla d\vert > 0 \quad 
\mbox{on $\ooo{\OOO}$}, \quad d=0 \quad\mbox{on $\ppp\OOO\setminus
\Gamma$}, \quad \nabla d\cdot \nu \le 0 \quad \mbox{on $\ppp\OOO\setminus 
\Gamma$}                         \eqno{(2.1)}
$$ 
(e.g, Imanuvilov \cite{Ima}).  Here $\nu$ denotes the unit outward normal 
vector to $\ppp\OOO$.

For arbitrarily fixed $t_0\in (0,T)$ and $\delta>0$ such that 
$0 < t_0 - \delta \le t_0 + \delta < T$, we set
$$
I = (t_0-\delta, t_0+\delta), \quad Q_I = \OOO \times I.
$$

We set 
$$
P_k v(x,t):= \ppp_tv + (-1)^k a(x,t)\Delta v + R(x,t,v),
\quad k=1,2,
$$
where $a \in C^{2}(\ooo{Q_I})$, $>0$ on $\ooo{Q_I}$, and 
$$
\vert R(x,t,v)\vert \le C_0(\vert v(x,t)\vert 
+ \vert\nabla v(x,t)\vert), \quad (x,t)\in Q_I.
                                     \eqno{(2.2)}
$$
Moreover, let
$$
\va(x,t) = e^{\la(d(x) - \beta (t-t_0)^2)},
$$
where $\la>0$ is a sufficiently large parameter and $\beta>0$ is
arbitrarily given.
Henceforth $C>0$ denote generic constants which independent of 
$s>0$, but depends on $\la, \beta, C_0$ in (2.2).
Then
\\
{\bf Lemma 1.}
{\it
There exist constants $s_0>0$ and $C>0$ such that 
$$
 \int_{Q_I} \left( \frac{1}{s}(\vert \ppp_tv\vert^2
+ \vert \Delta v\vert^2) + s \vert\nabla v\vert^2 
+ s^3\vert v\vert^2 \right)\weight dxdt
\le Cs^4\int_{Q_I} \vert P_kv\vert^2 \weight dxdt
+ C\BBB(v), \quad k=1,2 \eqno{(2.3)}
$$
for all $s > s_0$ and $v\in H^{2,1}(Q_I)$ satisfying 
$v\in H^1(\ppp\OOO\times I)$.
Here and henceforth we set 
\begin{align*}
& \BBB(v) := e^{Cs}\Vert v\Vert^2_{H^1(\Gamma\times I)}
+ s^3\int_{(\ppp\OOO\setminus \Gamma)\times I}
(\vert v\vert^2 + \vert \nabla_{x,t}v\vert^2) e^{2s} dSdt\\
+& s^2\int_{\OOO} (\vert v(x,t_0-\delta)\vert^2
+ \vert \nabla v(x,t_0-\delta)\vert^2
+ \vert v(x,t_0+\delta)\vert^2
+ \vert \nabla v(x,t_0+\delta)\vert^2) e^{2s\va(x,t_0-\delta)} dx.
\end{align*}
}
The proof of the lemma with $k=1$ is done similarly to Lemma 7.1 
(p.186) in Bellassoued and Yamamoto \cite{BY} or 
Theorem 3.2 in Yamamoto \cite{Y09} by keeping all the boundary 
integral terms $v\vert_{\ppp Q}$ which are produced by integration by parts 
and using $d\vert_{\ppp\OOO\setminus \Gamma}=0$ in (2.1).
The proof for $k=2$ follows directly from the case $k=1$ by setting 
$V(x,t):= v(x,2t_0-t)$ and using $\va(x,t) = \va(x,2t_0-t)$ for 
$(x,t) \in Q_I$.

We emphasize that the backward parabolic Carleman estimate is the same as
the forward parabolic Carleman estimate thanks to the symmetry of 
the weight $\va(x,t)$ with respect to $t$ centered at $t_0$.

Using the Carleman estimate (2.3) we  prove a Carleman estimate for a mean field game system.
Setting $y:= u - \wwwu$ and $z:= v-\wwwv$ and subtracting 
the system (1.1) with $(\wwwu,\wwwv,\www{F},\www{G})$ from
(1.1) with $(u,v,F,G)$, we reach
$$
\left\{ \begin{array}{rl}
& \ppp_ty + a_1(x,t)\Delta y + R_1(x,t,y) = h(x,t)z + F-\www{F},\\
& \ppp_tz - a_2(x,t)\Delta z + R_2(x,t,z) = \kappa v\Delta y + R_3(x,t,y)
+ G - \www{G}              \quad \mbox{in $Q_I$}.
\end{array}\right.
                 \eqno{(2.4)}
$$
Here by (1.2) and (1.3), we can verify 
$$
\vert R_j(x,t,y)\vert \le C_0\sum^1_{k=0} \vert \nabla^k y(x,t)\vert,
\quad j=1,3, \quad
\vert R_2(x,t,z)\vert \le C_0\sum^1_{k=0} \vert \nabla^k z(x,t)\vert,
\quad (x,t) \in Q_I.
                                   \eqno{(2.5)}
$$
We apply Carleman estimate (2.3) to  the first equation in (2.4), and multiply the resulting equality by $s$: 
$y$ and obtain
$$
 \int_{Q_I} ( \vert \ppp_ty\vert^2 + \vert \Delta y\vert^2 
+ s^2 \vert\nabla y\vert^2 + s^4\vert y\vert^2 )\weight dxdt
$$
$$
\le C\int_{Q_I} s\vert hz\vert^2 \weight dxdt + C\int_{Q_I} 
s\vert F-\www{F}\vert^2 \weight dxdt + Cs\BBB(y) 
                                                  \eqno{(2.6)}
$$
for all $s > s_0$.  In terms of (2.5), application of (2.3) with $k=1$ to 
$z$ yields
$$
 \int_{Q_I} \left( \frac{1}{s}(\vert \ppp_tz\vert^2 + \vert \Delta z\vert^2) 
+ s\vert\nabla z\vert^2 + s^3\vert z\vert^2 \right)\weight dxdt
$$
$$
\le C\int_{Q_I} (\vert \kappa\Delta y\vert^2 + \vert y\vert^2
+ \vert \nabla y\vert^2) \weight dxdt 
+ C\int_{Q_I} \vert G-\www{G}\vert^2 \weight dxdt + C\BBB(z) 
                                                  \eqno{(2.7)}
$$
for all $s > s_0$.

Using $\kappa \in L^{\infty}(Q_I)$ and substituting (2.6) into the terms 
including $\Delta y, \nabla y, y$ on the right-hand side of (2.7), we have
\begin{align*}
& \int_{Q_I} \left( \frac{1}{s}(\vert \ppp_tz\vert^2 + \vert \Delta z\vert^2) 
+ s\vert\nabla z\vert^2 + s^3\vert z\vert^2 \right)\weight dxdt\\
\le& C\int_{Q_I} s\vert z\vert^2 \weight dxdt 
+ C\int_{Q_I} (s\vert F-\www{F}\vert^2 + \vert G-\www{G}\vert^2) \weight dxdt 
+ Cs(\BBB(y) + \BBB(z))
\end{align*}
for all large $s>0$.
Absorbing the first term on the right-hand side into the left-hand side by 
choosing $s>0$ sufficiently large, we see
$$
 \int_{Q_I} \left( \frac{1}{s}(\vert \ppp_tz\vert^2 + \vert \Delta z\vert^2) 
+ s\vert\nabla z\vert^2 + s^3\vert z\vert^2 \right)\weight dxdt
$$
$$
\le C\int_{Q_I} (s\vert F-\www{F}\vert^2 + \vert G-\www{G}\vert^2)\weight dxdt 
+ Cs(\BBB(y) + \BBB(z))                                          \eqno{(2.8)}
$$
for all $s>s_0$.
Adding (2.8) and (2.6) and choosing $s>0$ large again to absorb the term
$\int_{Q_I} s\vert hz\vert^2 \weight dxdt$ on the 
right-hand side into the left-hand side, we obtain
\\
{\bf Theorem 2 (Carleman estimate for a mean field game)}.
{\it
There exist constants $s_0>0$ and $C>0$ such that 
\begin{align*}
& \int_{Q_I} \biggl( \vert \ppp_t(u-\www{u})\vert^2 
+ \vert (\Delta (u-\wwwu)\vert^2 
+ s^2\vert\nabla (u-\wwwu)\vert^2 + s^4\vert u-\wwwu\vert^2 
+ \frac{1}{s}(\vert \ppp_t(v-\www{v})\vert^2 
+ \vert \Delta (v-\wwwv)\vert^2) \\
+ &s\vert\nabla (v-\wwwv)\vert^2 + s^3\vert v-\wwwv\vert^2 \biggr)\weight dxdt
\le C\int_{Q_I} (s\vert F-\www{F}\vert^2 + \vert G-\www{G}\vert^2)\weight dxdt\\+ & Cs(\BBB(u-\wwwu) + \BBB(v-\wwwv)) \quad \mbox{for all $s>s_0$}.
\end{align*}
}
\section{Proof of Theorem 1}

We arbitrarily choose $t_0\in (0,T)$ and $\delta>0$ such that 
$0<t_0-\delta<t_0+\delta<T$.
We define 
$$
d_0:= \min_{x\in \ooo{\OOO}} d(x),\quad
d_1:= \max_{x\in \ooo{\OOO}} d(x),\quad
0 < r < \left( \frac{d_0}{d_1}\right)^{\hhalf} < 1.
                                         \eqno{(3.1)}
$$
We note that $0<r<1$.

We now show
\\
{\bf Lemma 2.}
{\it 
Under regularity condition (1.3), if $u=\wwwu$, $v=\wwwv$, 
$\nabla u = \nabla\wwwu$ and $\nabla v = \nabla \wwwv$ on
$\Gamma \times (t_0-\delta,\, t_0+\delta)$ imply
$u=\wwwu$ and $v=\wwwv$ in $\OOO \times (t_0-r\delta,\, t_0+r\delta)$.
}

For the proof of Theorem 1, it suffices to prove Lemma 2.
Indeed, since $t_0 \in (0,T)$ and $\delta>0$ can be arbitrarily 
chosen and the Carleman estimate is invariant 
with respect to $t_0$ provided that $0<t_0-\delta<t_0+\delta<T$, we can apply
Lemma 2 by changing $t_0$ over $(\delta, T-\delta)$ to obtain
$u=\wwwu$ and $v=\wwwv$ in $\OOO \times ((1-r)\delta,\, T-(1-r)\delta)$.
Since $\delta>0$ can be arbitrary, this means that 
$u=\wwwu$ and $v=\wwwv$ in $\OOO\times (0,T)$.
\\
{\bf Proof of Lemma 2.}
Once we derived the relevant Carleman estimate in Theorem 2, the proof of
Lemma 2 is done similarly to Proposition 2 in \cite{HIY} as follows.
First we determine the constant $\beta > 0$ in the weight of the Carleman 
estimate such that 
$$
\frac{d_1-d_0}{\delta^2 - r^2\delta^2} < \beta < \frac{d_0}{r^2\delta^2}.
                                      \eqno{(3.2)}
$$
Here we note that (3.1) verifies
$0 < \frac{d_1-d_0}{\delta^2 - r^2\delta^2} < \frac{d_0}{r^2\delta^2}$, which 
allows us to choose $\beta$ satisfying (3.2).
  
For short descriptions, we set
\begin{align*}
& M_1:= \sum_{k=0}^1 
(\Vert\nabla^k_{x,t}(u-\wwwu)\Vert^2
_{L^2((\ppp\OOO\setminus \Gamma)\times I)}
+ \Vert\nabla^k_{x,t}(v-\wwwv)\Vert^2
_{L^2((\ppp\OOO\setminus \Gamma)\times I)}), \\
& M_2:= \sum_{k=0}^1 (\Vert (u - \wwwu)(\cdot,t_0 + (-1)^k\delta) \Vert^2
_{H^1(\OOO)} 
+ \Vert (v - \wwwv)(\cdot,t_0 + (-1)^k\delta) \Vert^2_{H^1(\OOO)}
\end{align*}
and $\mu_1:= e^{\la(d_1-\beta\delta^2)}$.
Since $u=\wwwu$ and $v=\wwwv$ on $\Gamma \times I$, Theorem 2 yields
$$
s^3\int_{Q_I} (\vert u-\wwwu\vert^2 + \vert v-\wwwv\vert^2)
\weight dxdt
\le Cs^5M_1e^{2s} + Cs^5M_2e^{2s\mu_1}
$$
for all large $s>0$.
We shrink the integration region of the left-hand side to 
$\OOO \times (t_0-r\delta,\, t_0+r\delta)$.
Then, since 
$\va(x,t) = e^{\la(d(x) - \beta(t-t_0)^2)}
\ge e^{\la(d_0-\beta r^2\delta^2)}=:\mu_2$ in 
$\OOO \times (t_0-r\delta,\, t_0+r\delta)$,
we obtain
$$
e^{2s\mu_2}\int_{\OOO \times (t_0-r\delta,\, t_0+r\delta)}
(\vert u-\wwwu\vert^2 + \vert v-\wwwv\vert^2) dxdt 
\le Cs^2M_1e^{2s} + Cs^2M_2e^{2s\mu_1},     
$$
that is,
$$
\Vert u-\wwwu\Vert^2_{L^2(\OOO \times (t_0-r\delta,\, t_0+r\delta))}
+ \Vert v-\wwwv\Vert^2_{L^2(\OOO \times (t_0-r\delta,\, t_0+r\delta))}
\le Cs^2M_1e^{-2s(\mu_2-1)} + Cs^2M_2e^{-2s(\mu_2-\mu_1)}        \eqno{(3.3)}
$$
for all large $s>0$.
Here, by (3.2), we see that 
$\mu_2 > \max\{ 1, \, \mu_1\}$,
and so we let $s \to \infty$ in (3.3), so that 
$u=\wwwu$ and $v = \wwwv$ in $\OOO \times (t_0-r\delta,\, t_0+r\delta)$.
Thus the proof of Lemma 2, and so Theorem 1 are complete.
$\blacksquare$

{\bf Acknowledgments.}
The work was supported by Grant-in-Aid for Scientific Research (A) 20H00117 
of Japan Society for the Promotion of Science.

\end{document}